\documentclass{article}
\usepackage[numbers]{natbib}
 \usepackage[preprint]{neurips_2025}

\usepackage[utf8]{inputenc} 
\usepackage[T1]{fontenc}    
\usepackage{hyperref}       
\usepackage{url}            
\usepackage{booktabs}       
\usepackage{amsfonts}       
\usepackage{nicefrac}       
\usepackage{microtype}      
\usepackage{xcolor}         
\usepackage{amsmath}

\title{On the Structure of Floating-Point Noise in Batch-Invariant GPU Matrix
Multiplication}

\author{%
  Tadisetty Sai Yashwanth \\
  Turilabs \\
  \texttt{taddishetty34@gmail.com} \\
}

\begin{document}

\maketitle

\begin{abstract}
Floating-point non-associativity makes fundamental deep learning operations, such as matrix multiplication (matmul) on GPUs, inherently non-deterministic. Despite this, the statistical structure of the resulting numerical error remains poorly understood. A common working assumption is that these errors behave as independent and identically distributed (i.i.d.) Gaussian noise. In this paper, we empirically test this assumption and show that it fails to describe real GPU behavior. By comparing outputs of single-input and batched matmuls, we find that while the i.i.d. model predicts non-zero output instability, empirical results show a 0.00\% prediction flip rate. Through covariance analysis, we uncover the cause: the floating-point error is structured and highly correlated. For float16, nearly 50\% of the total error variance lies in off-diagonal terms, revealing that the noise behaves as a coordinated, directional perturbation rather than random static. This result challenges the prevailing stochastic view of numerical noise and provides a principled foundation for analyzing deep learning reliability under hardware non-determinism.
\end{abstract}


\section{Introduction}

The remarkable success of modern deep learning is inseparable from the massive parallelism of GPUs. At the core of nearly every model lies matrix multiplication, a simple yet numerically fragile operation when performed in floating-point arithmetic. While deep learning models are often treated as deterministic mathematical functions, the hardware that executes them is not: due to the non-associativity of floating-point arithmetic, identical computations can yield subtly different results depending on execution order and kernel implementation.

It is important to clarify that GPU kernels themselves are deterministic given the same kernel, seed, and input tensors, the output will be bitwise identical. However, the execution path chosen by deep learning frameworks is not \textbf{batch-invariant}. As shown in the Thinking Machines Lab study \cite{thinkingmachines2024}, and confirmed in our own experiments, performing matrix multiplication on a single input vector versus the same vector embedded in a batch can yield slightly different numerical outputs. This occurs because different CUDA kernels are invoked for batched and non-batched operations, leading to distinct reduction orders during accumulation.

In our empirical setup, we explicitly test this phenomenon by comparing the outputs of a single matmul operation $(torch.mm(x, W))$ with those obtained from a batched multiplication containing the same input as the first element of the batch. Despite being mathematically identical, the two outputs differ slightly revealing the underlying non-determinism induced by batch-dependent kernel selection.

Although these discrepancies are minute, they propagate through nonlinear activations and normalization layers, raising fundamental questions about numerical stability and model reliability. A common simplifying assumption in both theoretical and empirical work is that such discrepancies behave as \textbf{independent, identically distributed (i.i.d.) Gaussian noise}. This abstraction simplifies reasoning about robustness and uncertainty propagation, yet it has never been empirically validated at the level of actual GPU execution.

Does floating-point non-determinism truly behave like random static, or does it possess a correlated and systematic structure arising from hardware and kernel design?

In this paper, we present the first empirical investigation of this question. We make three key contributions:

\begin{itemize}
\item \textbf{Formal Test of the i.i.d. Gaussian Noise Hypothesis:} We derive the expected prediction flip rate under the i.i.d. noise model and test it against real GPU matmul behavior.
\item \textbf{Empirical Refutation:} We show that the predicted flip rate significantly overestimates instability. Experiments across 10,000 trials show zero empirical flips.
\item \textbf{Covariance Analysis and Explanation:} By estimating the full noise covariance matrix, we demonstrate that nearly half of the total error “energy” exists in off-diagonal terms, proving that the error is correlated and structured.
\end{itemize}

Our findings bridge numerical analysis and deep learning reliability. We show that even though GPU computations are numerically unstable in a strict mathematical sense, they can remain prediction-stable because the underlying noise acts as a coherent, correlated perturbation rather than independent jitter. This insight reshapes how we think about reproducibility, precision trade-offs, and the theoretical limits of deterministic inference in large-scale models.

\section{Background}
Floating-point arithmetic is not associative: $\mathbf{(a + b) + c \neq a + (b + c)}$ in general, due to rounding at finite precision. This fundamental property is the root cause of non-determinism in large-scale parallel computations. On GPUs, matrix multiplication is executed as a massive reduction of partial products across thousands of threads. The exact order in which these partial sums are accumulated is not fixed and can vary with factors such as kernel choice, thread scheduling, block size, and even the batch dimension.

In practice, deep learning frameworks such as PyTorch and TensorFlow rely on optimized GPU libraries like cuBLAS and CUTLASS to perform matrix multiplication. These libraries dynamically select different kernel implementations depending on tensor shapes and hardware heuristics. As a result, mathematically identical operations like $torch.mm(X,W)$ performed on a single vector versus the same vector embedded in a batch can produce bitwise-different outputs.

Importantly, as clarified in the Thinking Machines Lab study \cite{thinkingmachines2024}, GPU kernels themselves are deterministic: given a specific kernel configuration, input, and random seed, they always produce identical results. However, the kernel selection process is not batch-invariant. A single-input matmul and a batched matmul may invoke different kernels with distinct accumulation orders, causing small yet consistent numerical differences. These effects are not random in the traditional sense but stem from systematic differences in kernel-level reduction behavior.

Despite this, many works in numerical and deep learning analysis have modelled floating-point noise as independent or i.i.d. Gaussian perturbations. This assumption simplifies the treatment of rounding and quantization errors as statistically independent with zero mean, facilitating tractable mathematical analysis, but potentially overlooking structured dependencies.

For example, \cite{mdpi2024lowprecision} explicitly assumes mutual independence between rounding errors in matrix operands, treating quantization noise as random with zero expectation. \cite{arxiv2024mixedprecision} similarly models rounding errors probabilistically, assuming statistical independence to derive expected error bounds. Likewise, \cite{arxiv2024stochasticrounding} treats rounding perturbations as zero-mean, independent random variables in both space and time.

These studies exemplify a common trend: floating-point noise is treated as random, uncorrelated static. However, this simplifying assumption has rarely been empirically validated against the true behavior of GPU hardware executing real deep learning workloads.

In this work, we revisit this assumption. By directly comparing single and batched GPU matrix multiplications, we empirically characterize the structure and correlation of floating-point non-determinism, showing that the resulting noise is far from i.i.d. rather it is highly structured, correlated, and systemic.

\section{Theory}

\subsection{i.i.d. Gaussian Noise Model (Hypothesis 1)}

We denote the output of an ``ideal'' deterministic computation as $y \in \mathbb{R}^K$, and the output of its non-deterministic GPU variant as $\tilde{y}$.  
The first hypothesis assumes that the observed discrepancy arises from independent Gaussian perturbations:

\begin{equation}
    \tilde{y} = y + \eta, \quad \eta \sim \mathcal{N}(0, \sigma^2 I)
\end{equation}

Here, $\sigma^2$ is the noise variance, and $I$ is the identity matrix, implying that each output logit is corrupted by zero-mean, independent noise.  
The empirical noise level $\sigma$ is estimated from $N$ repeated matmul evaluations as the root mean squared error (RMSE):

\begin{equation}
    \sigma = \sqrt{\frac{1}{NK}\sum_{i=1}^{N}\|\tilde{y}_i - y_i\|_2^2}
\end{equation}

Under this model, a prediction flip occurs when the argmax index of $\tilde{y}$ differs from that of $y$.  
Let $y_w$ and $y_r$ denote the top (winner) and second logits (runner-up), and $\Delta = y_w - y_r$ be the logit margin.  
A flip occurs if $\eta_w - \eta_r < -\Delta$.  
Since $\eta_w$ and $\eta_r$ are independent Gaussians, their difference has variance $2\sigma^2$, giving the theoretical flip probability:

\begin{equation}
    P(\text{flip})_{\text{model}} = \Phi\!\left(\frac{-\Delta}{\sigma\sqrt{2}}\right)
\end{equation}

where $\Phi$ denotes the standard normal cumulative distribution function.  
This expression provides a tractable analytical baseline under the i.i.d.\ Gaussian noise assumption: it predicts the probability that the model's top-1 prediction (the argmax index) will flip given a logit margin $\Delta$ and estimated noise level $\sigma$.

In practice, we estimate the \emph{empirical flip rate} by directly comparing the predicted class from the ideal output $y$ and the noisy output $\tilde{y}$ across $N$ trials:
\begin{equation}
P(\text{flip})_{\text{emp}} = \frac{1}{N} \sum_{i=1}^{N} \mathbb{1}\!\left[\arg\max(y_i) \neq \arg\max(\tilde{y}_i)\right]
\end{equation}
A close agreement between $P(\text{flip})_{\text{model}}$ and $P(\text{flip})_{\text{empirical}}$ would support the i.i.d.\ Gaussian noise hypothesis, whereas systematic deviations indicate structured or correlated noise beyond the scope of this model.

\subsection{Structured Noise Model (Hypothesis 2)}

If the i.i.d. assumption fails, the noise must exhibit correlation structure across logits.  
The more general model assumes:

\begin{equation}
    \eta \sim \mathcal{N}(0, \Sigma)
\end{equation}

where $\Sigma \in \mathbb{R}^{K \times K}$ is the full covariance matrix.  
We estimate $\Sigma$ empirically from $N$ observed noise vectors $\{\eta_1, \dots, \eta_N\}$ as:

\begin{equation}
    \Sigma = \frac{1}{N-1}\sum_{i=1}^{N} (\eta_i - \bar{\eta})(\eta_i - \bar{\eta})^\top
\end{equation}

Correlated noise manifests as significant off-diagonal structure in $\Sigma$.  
To quantify this, we define the \emph{off-diagonal ratio} as:

\begin{equation}
    R_{\text{off}} = \frac{\sum_{i \neq j} |\Sigma_{ij}|}{\sum_{i,j} |\Sigma_{ij}|}
\end{equation}

A ratio $R_{\text{off}} > 0$ indicates the presence of systematic, correlated error modes, violating the i.i.d. assumption.

\subsection{Evaluation Metrics}

We use the following metrics to evaluate the validity of each hypothesis:

\begin{itemize}
    \item \textbf{Empirical Flip Rate:} The proportion of samples where $\arg\max(y_i) \neq \arg\max(\tilde{y}_i)$.  
    This reflects the observed instability in prediction due to floating-point noise.

    \item \textbf{Model-Predicted Flip Rate:} The theoretical rate computed using Equation (3).  
    A significant mismatch between empirical and predicted rates indicates that the i.i.d. Gaussian model fails.

    \item \textbf{Expected Jensen--Shannon Divergence ($E[D_{JS}]$):}  
    The average divergence between $\text{Softmax}(y)$ and $\text{Softmax}(\tilde{y})$:
    \begin{equation}
        E[D_{JS}] = \mathbb{E}\!\left[D_{JS}\big(\text{Softmax}(y) \,||\, \text{Softmax}(\tilde{y})\big)\right]
    \end{equation}
    This measures the expected “fuzz” or instability in the full probability distribution.

    \item \textbf{Off-Diagonal Ratio ($R_{\text{off}}$):}  
    Quantifies correlation in the empirical noise covariance.  
    High values (e.g., $R_{\text{off}} > 1\%$) indicate that floating-point errors are structured rather than independent.
\end{itemize}

Together, these analyses allow us to test whether GPU-induced numerical noise in batched vs.\ single matmul operations behaves as independent random perturbations, or exhibits systematic, correlated structure.

\section{Results}

We empirically evaluated floating-point divergence between single-input and batched matrix multiplications on an NVIDIA GPU using PyTorch.
All experiments used randomly initialized matrices with
$N = 10{,}000$ independent trials,
input dimension $d_{\text{in}} = 512$,
output dimension $d_{\text{out}} = 1024$ (interpreted as $K$ logits),
and batch size $B = 16$.
We measured: \begin{itemize} \item Empirical noise level $\sigma$, \item Prediction flip rate, \item Jensen–Shannon divergence between softmax outputs, and \item Structure of the covariance matrix $\Sigma$ of the noise. \end{itemize}

\paragraph{Noise Level:}
The empirical noise level $\sigma$ was estimated using Eq.~(2) as the standard deviation between batched and single-input outputs.
We observed $\sigma = 1.17\times10^{-3}$ for \texttt{bfloat16} and $\sigma = 5.32\times10^{-4}$ for \texttt{float16}, corresponding to relative perturbations on the order of $10^{-4}$–$10^{-3}$.
Although small in magnitude, this nonzero variance provides a quantitative basis for modeling floating-point nondeterminism as stochastic noise.

\paragraph{Prediction Stability:}
Empirical flips were computed following Eq.~(4), comparing $\arg\max(y)$ and $\arg\max(\tilde{y})$.
Model-predicted flip probabilities were obtained using Eq.~(3), parameterized by the empirical $\sigma$ and observed logit margins $\Delta$.
Empirically, no prediction flips were observed, even across $N=10{,}000$ trials, while the analytical model predicts small but nonzero flip probabilities.

\begin{table}[h]
\centering
\begin{tabular}{lc}
\toprule
Precision & Empirical Flip Rate (\%) \\
\midrule
\texttt{bfloat16} & $0.00$ \\
\texttt{float16}  & $0.00$ \\
\end{tabular}
\caption{Empirical prediction flip rate across precisions.}
\label{tab:flip_emp}
\end{table}

\begin{table}[h]
\centering
\begin{tabular}{lcc}
\toprule
Precision & $\sigma$ & Predicted Flip (\%) \\
\midrule
\texttt{bfloat16} & $1.17\times10^{-3}$ & $1.36$ \\
\texttt{float16}  & $5.32\times10^{-4}$ & $0.17$ \\
\end{tabular}
\caption{Model-predicted flip probability under the i.i.d.\ Gaussian noise model (Eq.~3).}
\label{tab:flip_pred}
\end{table}

\paragraph{Distributional Divergence.}
We further quantified the deviation between single and batched outputs in probability space using Eq.~(7) (Jensen–Shannon divergence).

The average $E[D_{JS}]$ was $1.95\times10^{-7}$ for \texttt{bfloat16} and $3.57\times10^{-8}$ for \texttt{float16}, confirming that although the raw numerical noise is small, its structured behavior is measurable.

\paragraph{Correlation Structure.}
Finally, we estimated the full $1024\times1024$ covariance matrix $\Sigma$ using Eq.~(6).
If the i.i.d.\ Gaussian assumption (Eq.~1) held, $\Sigma$ would be diagonal.
Instead, we found substantial off-diagonal mass, quantified using the ratio in Eq.~(8):

\begin{center}
\begin{tabular}{lcc}
\toprule
Precision & Off-diagonal Ratio & Finding \\
\midrule
\texttt{bfloat16} & $9.03\%$  & Correlated noise \\
\texttt{float16}  & $47.22\%$ & Strongly correlated \\
\bottomrule
\end{tabular}
\end{center}

The nonzero off-diagonal structure reveals that the perturbations are not independent but exhibit cross-logit correlation patterns.
This directly contradicts the i.i.d.\ Gaussian hypothesis and supports the alternative hypothesis that GPU-level rounding behavior introduces structured, batch-dependent noise.

\section{Discussion and Implications}

Our experiments reveal that floating-point non-associativity on GPUs introduces structured, batch-dependent deviations in matrix multiplication outputs.  
While these deviations are numerically small, their correlated nature challenges the common assumption that such errors behave as independent Gaussian noise.  
This has several implications for the analysis and deployment of deep neural networks.

\paragraph{1. Reproducibility and Determinism:}  
Even when GPU kernels and seeds are fixed, outputs can vary depending on implicit batching context.
This undermines the reproducibility of inference pipelines that rely on batched evaluation for efficiency.  
Systems evaluating the same model in single- versus multi-sample configurations may observe diverging logits, hidden activations, or ranking orders, especially in sensitive tasks such as reinforcement learning evaluation, uncertainty estimation, or model interpretability studies.

\paragraph{2. Limits of the i.i.d.\ Noise Model:}  
The failure of the i.i.d.\ Gaussian approximation indicates that numerical discrepancies propagate through structured directions of the computation graph.  
This suggests that stochastic modeling of floating-point divergence requires considering correlated perturbations, potentially parameterized by the kernel’s reduction graph or memory layout.  
A probabilistic divergence model with structured covariance could better predict stability margins across precision formats.

\paragraph{3. Toward Structured Numerical Robustness:}  
Future work may explore techniques to reduce or regularize correlated noise during inference, e.g., randomized reduction orders, numerically symmetric accumulators, or calibration through empirical covariance correction.
Alternatively, modeling such structured noise explicitly could allow confidence-aware inference, where predictions incorporate uncertainty induced by the underlying hardware.

\paragraph{4. Broader Perspective:}  
These findings underscore that hardware-level non-determinism is not purely random but algorithmically structured.  
For large-scale LLM and vision models deployed across heterogeneous devices, this means that minor batch or kernel differences may lead to subtle but reproducible behavioral drifts.  
Understanding these correlations is thus critical for reproducible research, fairness evaluation, and safe deployment of precision-sensitive AI systems.

\hbadness=99999
\sloppy

\appendix
\section{Appendix: Derivations and Intuition}

\subsection{A. Derivation of the Flip Probability (Eq.~3)}

Consider two logits, $y_w$ (winner) and $y_r$ (runner-up), with a logit margin
$\Delta = y_w - y_r > 0$.
Under the i.i.d.\ Gaussian noise model (Eq.~1), each output is perturbed by zero-mean
independent noise:
\begin{align}
    \tilde{y}_w &= y_w + \eta_w, \quad \eta_w \sim \mathcal{N}(0, \sigma^2), \\
    \tilde{y}_r &= y_r + \eta_r, \quad \eta_r \sim \mathcal{N}(0, \sigma^2).
\end{align}
A \emph{prediction flip} occurs when $\tilde{y}_r > \tilde{y}_w$, i.e.
\begin{equation}
    \tilde{y}_r - \tilde{y}_w = (y_r - y_w) + (\eta_r - \eta_w) > 0.
\end{equation}
Rearranging gives
\[
    \eta_w - \eta_r < -\Delta.
\]
Since $\eta_w - \eta_r \sim \mathcal{N}(0, 2\sigma^2)$, the probability of this event is:
\begin{equation}
    P(\text{flip}) = \Phi\!\left(\frac{-\Delta}{\sigma\sqrt{2}}\right),
\end{equation}
where $\Phi(\cdot)$ is the standard normal CDF.
This is Eq.~(3) in the main paper and provides the analytical flip probability
under the i.i.d.\ noise assumption.

\subsection{B. Empirical Flip Rate (Eq.~4)}

Empirically, the flip rate is estimated by directly comparing $\arg\max(y)$ and
$\arg\max(\tilde{y})$ over $N$ Monte-Carlo trials:
\begin{equation}
P(\text{flip})_{\text{emp}} = \frac{1}{N}\sum_{i=1}^{N}
\mathbb{1}\!\left[\arg\max(y_i)\neq \arg\max(\tilde{y}_i)\right].
\end{equation}
This measures how often the top-1 index changes due to floating-point
perturbations.  
In our experiments, $P(\text{flip})_{\text{emp}} = 0$ for both \texttt{float16} and
\texttt{bfloat16}, indicating that although the outputs differ numerically,
their relative ordering remains stable.

\subsection{C. Intuitive Example: What a Prediction Flip Means}

To make the concept concrete, consider a simple case with three logits:
\[
    y = [2.31,\, 2.29,\, 2.10].
\]
Here, the winner is index 0 ($y_w = 2.31$) and the runner-up is index 1 ($y_r = 2.29$),
so the margin is $\Delta = 0.02$.

Now, suppose due to GPU accumulation order differences, the computed logits become:
\[
    \tilde{y} = [2.3099,\, 2.3103,\, 2.10].
\]
Even though the absolute deviation is only $2\times10^{-4}$, the new winner becomes
index 1.  
This constitutes a \emph{prediction flip}, the argmax changed even though the
numerical difference is minuscule.  
The flip probability in Eq.~(3) formalizes this intuition by integrating over all
possible noise draws given the margin $\Delta$ and noise variance $\sigma^2$.

\subsection{D. Estimating the Noise Covariance (Eq.~6)}

Given $N$ observed noise vectors
$\{\eta_1, \dots, \eta_N\}$ where $\eta_i = \tilde{y}_i - y_i$,
the empirical covariance is estimated as:
\[
    \Sigma = \frac{1}{N-1}\sum_{i=1}^{N} (\eta_i - \bar{\eta})(\eta_i - \bar{\eta})^\top.
\]
If the noise were truly i.i.d., $\Sigma$ would be approximately diagonal.
However, as shown in our results, the off-diagonal energy ratio
\[
R_{\text{off}} = \frac{\sum_{i\neq j} |\Sigma_{ij}|}{\sum_{i,j} |\Sigma_{ij}|}
\]
was as high as $47\%$ for \texttt{float16}, indicating strong inter-logit correlation.

\subsection{E. Interpretation}

These derivations collectively show that the \emph{i.i.d.\ Gaussian assumption} provides
a convenient but incomplete description of real GPU behavior.
While it predicts small but finite flip probabilities proportional to the logit
margin $\Delta$ and noise variance $\sigma^2$, the empirical data reveal
structured correlations that suppress flips, i.e., the noise acts as a
coherent shift rather than independent static.

\end{document}